\documentclass{amsproc}

\usepackage{graphicx}
\usepackage{cite}
\usepackage{url}
\newtheorem{theorem}{Theorem}[section]
\newtheorem{lemma}[theorem]{Lemma}

\numberwithin{equation}{section}

\newcommand{\abs}[1]{\lvert#1\rvert}

\begin{document}

\title{Separating Thickness from Geometric Thickness}

\author{David Eppstein}
\address{School of Information and Computer Science\\
University of California, Irvine\\
Irvine, CA 92697-3425, USA}
\email{eppstein@ics.uci.edu}



\thanks{Work supported in part by NSF grant CCR-9912338. A preliminary version of this work appeared in the Proceedings of the 10th International Symposium on Graph Drawing (GD '02), Springer-Verlag Lecture Notes in Computer Science 2528, 2002, pp. 150--162.}

\thanks{This paper is in final form and no version of it will be submitted for publication elsewhere.}

\begin{abstract}
We show that graph-theoretic thickness
and geometric thickness are not
asymptotically equivalent: for every~$t$, there exists a graph with
thickness three and geometric thickness at least $t$.
\end{abstract}

\maketitle

\section{Introduction and Main Result}

The graph-theoretic property of {\em thickness} and its variants formalize the concept of {\em layered graph drawing}, in which a graph must be drawn in multiple layers (or with multiple colors) so that each edge belongs to a single layer and no two edges in the same layer cross.  The vertices of the graph should be represented as points and should exist in the same positions in each layer.  There are three important variants of thickness:

\begin{enumerate}
\item
The {\em thickness} of a graph $G$, denoted $\theta(G)$, is the minimum number of layers for which a drawing of $G$ exists, without restriction on the number of bends per edge~\cite{Kai-AMSUH-73}.

\item
The {\em geometric thickness} of a graph $G$, denoted $\bar\theta(G)$, also known as {\em real linear thickness} or {\em rectilinear thickness}, is the minimum number of layers for which a drawing exists without any bends in the edges: each edge must be represented as a straight line segment~\cite{DilEppHir-JGAA-00,Kai-AMSUH-73}.

\item 
The {\em book thickness} of a graph $G$, denoted ${\rm bt}(G)$,  is the minimum number of layers for which a drawing exists without any bends in the edges and with all vertices placed in convex position~\cite{BerKai-JCTB-79}.
\end{enumerate}

Since each variant imposes additional constraints on the previous one, these parameters are always ordered: $\theta(G)\le\bar\theta(G)\le {\rm bt}(G)$, and it is of interest to determine how tight these inequalities are.
Due to F\'ary's theorem~\cite{Far-ASMS-48}, the graphs with thickness one and geometric thickness one coincide (both are just the planar graphs), and are a strict superset of the graphs with book thickness one (the outerplanar graphs).  Thickness and geometric thickness were known to differ for infinitely many graphs: in particular, complete graphs on $n$ vertices have thickness $n/6+O(1)$~\cite{AleGon-MUS-76,Bei-GTTP-67,BeiHar-CJM-65,May-JCTB-72,Vas-ms-??} but have geometric thickness at least $n/5.646$~\cite{DilEppHir-JGAA-00}.

In previous work~\cite{math-CO-0109195} we showed that book thickness can not be bounded by any function of the geometric thickness: for any $t$ there exists a graph with geometric thickness two and book thickness at least $t$.  The proof idea is to consider the graphs
 $G_2(n)$ with geometric thickness two, having as vertices the singleton and doubleton subsets of an $n$-element set, with an edge between two subsets when one is contained in the other.  Ramsey's theorem shows that, if these graphs are partitioned into $k$ layers, for $k$ fixed and $n$ large, then some two layers include a $G_2(5)$.  But $G_2(5)$ is nonplanar and therefore has book thickness greater than two, so these two layers cannot be part of a book layout.  The same bound on the book thickness of $G_2(n)$, with the same proof, was independently discovered by Blankenship and Oporowski~\cite{BlaOpo-TR-99}.
 
In this paper, we prove a similar separation between thickness and geometric thickness:  for any $t$ there exists a graph with thickness three and geometric thickness at least $t$.
This implies that known strategies for approximating thickness (layered drawing with bends) will not lead to good approximation algorithms for layered drawing without bends.

The overall outline of our proof is very similar to the previous separation of geometric thickness from book thickness: we use subset inclusions to define a family of graphs with thickness three, show that
some graph in the family does not have geometric thickness three, and use Ramsey's theorem to conclude that the geometric thickness of the family is unbounded.
The difficult part of the proof is to show that some graph in the family does not have geometric thickness three, for which we use an intricate case analysis that also applies Ramsey theory multiple times in order to restrict the argument to drawings of certain uniform types.
Our main result is the following theorem.

\begin{theorem}
\label{thm:sepgeom}
For every $t$, a graph exists with
thickness three and geometric thickness at least $t$.
\end{theorem}

\begin{figure}[t]
\centering
\includegraphics[width=2.5in]{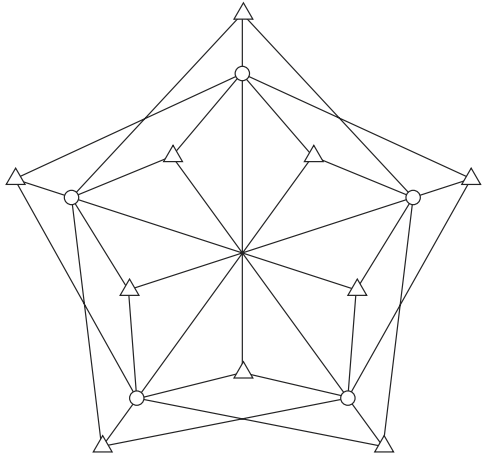}
\caption{Graph $G_3(5)$. The circular vertices represent singleton subsets, and the triangular vertices represent tripleton subsets.}
\label{fig:g35}
\end{figure}

\begin{proof}
For any $n>0$, define a graph $G_3(n)$ having as its $n+\binom{n}{3}$
vertices the singleton and tripleton subsets of an $n$-element set, with an edge between two subsets when one is contained in the other.  $G_3(4)$ forms the vertices and edges of a cube; $G_3(5)$ is depicted in Figure~\ref{fig:g35}. If we assign the edges of $G_3(n)$ to three layers, in such a way that the three edges incident to each tripleton are assigned to different layers, then each layer will form a forest of stars, which is a planar graph, so $G_3(n)$ has thickness at most three.  We show that there exists an $n_1$ such that $G_3(n_1)$ has geometric thickness more than three (Lemma~\ref{lem:notgt3}).

By Ramsey's theorem, there exists an $n_2$ such that, if the edges of $G_2(n_2)$ are given $t-1$ colors, one can find a three-colored $G_3(n_1)$ subgraph.
If the $(t-1)$-coloring is determined by a geometric embedding with $t-1$ layers, the subgraph would have a three-layer geometric embedding.  Since no such three-layer embedding exists, $G_3(n_2)$ does not have a $(t-1)$-layer geometric embedding.
\end{proof}

In the rest of this paper we fill in the details of the proof outline above.

\section{Restrictions to Coherent Drawings}

We use the following classical results of Ramsey theory~\cite{ErdSze-CM-35,GraRotSpe-90}.

\begin{lemma}[Ramsey's theorem]
\label{lem:Ram}
For every three positive integers $c$, $e$, and $\ell$ there is an integer
$R_e(\ell;c)$ with the following property:
if $\abs{S}\ge R_e(\ell;c)$ and the $e$-element subsets of $S$ are partitioned into $c$ classes, then there exists a set $T\subset S$
with $\abs{T}\ge\ell$ so that all $e$-element subsets of $T$ belong to the same class.
\end{lemma}

For instance, in the proof  of Theorem~\ref{thm:sepgeom} we choose
$c=\binom{t-1}{3}+\binom{t-1}{2}+t-1$,
$e=3$, and $\ell=n_1$.  The set $S$ is the set of singletons in $G_3(n_2)$, and the tripletons in $G_3(n_2)$ are partitioned into $c$ classes by the colors of their three incident edges.
Ramsey's theorem tells us that we can find a set $T$, $\abs{T}=n_1$, so that the subgraph $G_3(n_1)$ formed by the singleton and tripleton  subsets of $T$ uses the same three colors for each tripleton;
that is, we can find a three-colored $G_3(n_1)$ subgraph of $G_3(n_2)$.

\begin{lemma}[the Erd{\H{o}}s-Szekeres theorem]
\label{lem:Convex}
For every integer $k$ there is an integer $N(k)$ with the following property:
if $N(k)$ points are placed in the plane with no three points collinear, then
some subset of $k$ points forms the set of vertices of a convex polygon.
\end{lemma}

Define a {\em convex drawing} of $G_3(n)$ to be a drawing with geometric thickness three in which all singletons are in convex position.  We do not restrict the positions of the tripletons.

\begin{lemma}
\label{lem:condraw}
If every graph $G_3(n)$ has geometric thickness three, then every $G_3(n)$ has a convex drawing.
\end{lemma}

\begin{proof}
Sufficiently small perturbations of the vertex locations in a layered drawing will not create or remove crossings, so we can assume no three vertices are collinear.
To form a convex drawing of $G_3(n)$,
apply the Erd{\H{o}}s-Szekeres theorem to the singletons of $G_3(N(n))$, and form a
$G_3(n)$ subgraph from the resulting convex set of $n$ singletons and their adjacent tripletons.
\end{proof}

We number the singleton vertices of a convex drawing in clockwise order from $0$ to $n-1$, starting from an arbitrarily chosen vertex.  Using this numbering, we can partition the edges of $G_3(n)$ into three classes: a {\em low} edge connects a tripleton to the adjacent singleton with the smallest number, a {\em high} edge connects a tripleton to the adjacent singleton with the largest number, and a {\em middle} edge connects a tripleton to the remaining adjacent singleton.  Define a {\em coherent drawing} of $G_3(n)$ to be a convex drawing in which these three classes form the three layers of the drawing.

\begin{lemma}
\label{lem:coherent}
If every graph $G_3(n)$ has geometric thickness three, then every $G_3(n)$ has a coherent drawing.
\end{lemma}

\begin{proof}
By Lemma~\ref{lem:condraw} we can assume $G_3(n)$ has a convex drawing.
Partition the tripletons of such a drawing into 27 classes according to the layers to which the incident low, middle, and high edges belong.  By Lemma~\ref{lem:Ram} with $c=27$ and $e=3$, we can find a drawing in which each of the three sets of edges belongs to a single layer.
If this results in more than one set belonging to the same layer, we can move sets to distinct layers without introducing any crossings. 
\end{proof}

In a coherent drawing, we will denote tripletons by triples of symbols $abc$,
where each letter denotes the position of an adjacent singleton in the clockwise ordering of the drawing, and $a<b<c$ according to the clockwise ordering.

\section{Restriction to Outer Drawings}

Define an {\em inner drawing} of $G_3(n)$ to be a coherent drawing such that there exists a strictly convex curve $C$ on which all singletons lie, with all tripletons interior to $C$.  Similarly, define an {\em outer drawing} to be a coherent drawing such that there exists a strictly convex curve $C$ on which all singletons lie, with all tripletons exterior to $C$.  A drawing can be both inner and outer, for instance if all singletons and tripletons are in convex position.  Note that, in an inner drawing, all edges are strictly interior to $C$; however an edge in an outer drawing may have portions inside $C$ and portions outside $C$.

\begin{lemma}
\label{lem:onecross}
Suppose two tripletons $abc$ and $def$ are part of an inner drawing,
with $a<d<c<f$.  Then either the low edge of $abc$ crosses the high edge of $def$,
or the high edge of $abc$ crosses the low edge of $def$. However, both types of crossing can not occur simultaneously.
\end{lemma}

\begin{proof}
Consider the two paths formed by the low and high edges of each tripleton.
These paths cross $C$ and each path separates the endpoints of the other path; therefore the paths must cross an odd number of times.  However, each pair of edges can cross only once, and the only types of crossing permitted are those described in the lemma.  Therefore exactly one of these crossings occurs.
\end{proof}

We say that the two tripletons {\em cross convexly} if the low edge of $abc$ crosses the high edge of $def$, and that they {\em cross concavely} in the other case of Lemma~\ref{lem:onecross}.

\begin{figure}[t]
\centering
\includegraphics[width=4.5in]{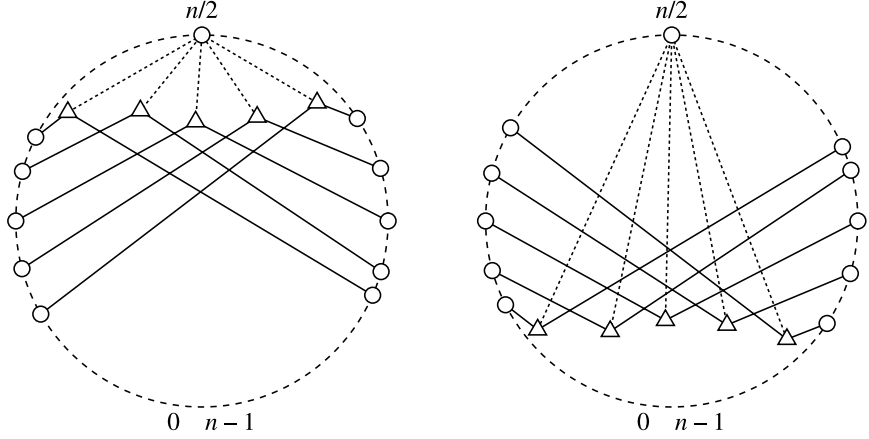}
\caption{Five tripletons crossing each other convexly (left) or concavely (right), from Lemma~\ref{lem:noinner}.}
\label{fig:innergrid}
\end{figure}

\begin{figure}[t]
\centering
\includegraphics[width=4.5in]{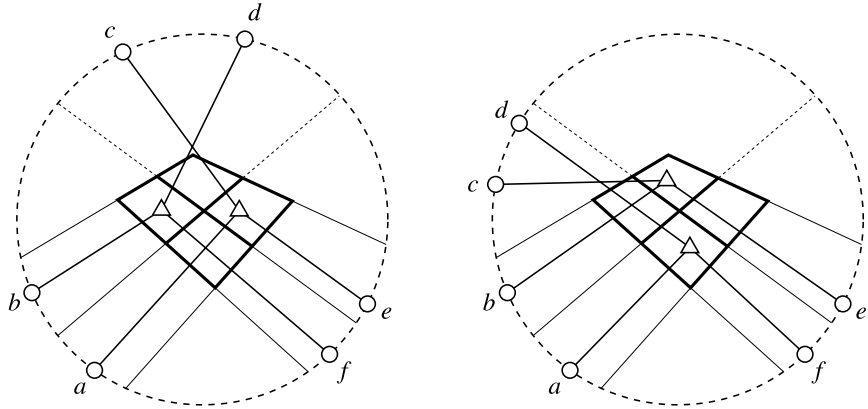}
\caption{Forced crossings in an inner drawing, from Lemma~\ref{lem:noinner}.}
\label{fig:innercross}
\end{figure}

\begin{lemma}
\label{lem:noinner}
There exists an $n_3$ such that $G_3(n_3)$ does not have an inner drawing.
\end{lemma}

\begin{proof}
Let $n$ be sufficiently large, and assume for convenience that $n=2$ (mod~$4$).  Consider the sequence of
tripletons ${2i,n/2,n/2+2i+2}$ for $0\le i<n/2-2$; each pair of tripletons meets the conditions of Lemma~\ref{lem:onecross}.  By applying Lemma~\ref{lem:Ram} with $c=2$, $e=2$, and $\ell=5$, for sufficiently large $n$, we can obtain a set of five such triples that all cross convexly or all cross concavely.  As Figure~\ref{fig:innergrid} shows, in either case, the low and high edges of these crossing triples form a grid graph, containing a complete $2\times 2$ subgrid.  (The quadrilaterals of this grid need not always be convex, but this does not affect our argument.)

Because we chose tripletons using only every other singleton from the ordering, it is now possible to find singletons $a$ and $b$ that lie between the three low edges forming the $2\times 2$ subgrid, and $e$ and $f$ that lie between the three high edges forming the subgrid (Figure~\ref{fig:innercross}).  Therefore, any tripleton with a low edge to $a$ or $b$ and a high edge to $e$ or $f$ must be forced to lie within a particular square of the subgrid.

Extend the two middle line segments of the subgrid into lines across curve $C$ (the thin dotted lines of the figure), partitioning $C$ into three curves: the curve clockwise of both extended lines, the curve counterclockwise of the extended lines, and the curve between the two lines.  The four points $n/2-2$, $n/2-1$, $n/2$, $n/2+1$ lie between $b$ and $e$, and by the pigeonhole principle we can find two of these points $c<d$ that are both within the same one of these three curves.

If $c$ and $d$ are both in the center curve of the three curves (shown in the left side of the figure) then tripletons $ace$ and $bdf$ must have crossing middle edges.
If $c$ and $d$ are both in the counterclockwise curve (shown on the right side of the figure), then tripletons $adf$ and $bce$ must have crossing middle edges; the case in which $c$ and $d$ belong to the clockwise curve is symmetric.  Thus, in all cases two tripletons have crossing middle edges, contradicting the assumption that we have a three-layer drawing with no crossings.
\end{proof}

\begin{lemma}
\label{lem:outer}
If every graph $G_3(n)$ has geometric thickness three, then every $G_3(n)$ has an outer  drawing.
\end{lemma}

\begin{proof}
Assume without loss of generality that $n$ is sufficiently large that, by Lemma~\ref{lem:noinner}, no inner drawing of $G_3(n)$ exists.
Let $r=R_3(n;2)$, as shown to exist in Lemma~\ref{lem:Ram}.
By Lemma~\ref{lem:coherent}, we can assume $G_3(r)$ has a coherent drawing.
Let $C$ be the convex hull of the singletons of this drawing, and perturb $C$ if necessary so that it remains convex and does not pass through any tripletons.
Partition the tripletons of the drawing into two classes according to whether they are inside or outside $C$, and apply Lemma~\ref{lem:Ram} with $c=2$ and $e=3$ to find a drawing of $G_3(n)$ in which all triples are inside $C$ or all triples are outside $C$.
Since no inner drawing exists, all triples are outside $C$, and we have an outer drawing of $G_3(n)$.
\end{proof}

\section{Classification of Outer Drawings}

In an outer drawing, the three angles formed by the three edges incident to a tripleton must include an angle greater than 180$^\circ$, and we can order the three edges clockwise starting at this large angle.  Define the {\em type} of a tripleton to be a symbol $xyz$, where $x$, $y$, and $z$ form a permutation of $0$, $1$, and $2$, and where the positions of the low, middle, and high number in the symbol are the same as the positions of the low, middle, and high edge in the clockwise ordering at the tripleton.  There are six possible types: $012$, $021$, $102$, $120$, $ 201$, and $210$.  If all tripletons in a drawing have the same type $xyz$, we say that the drawing is of type $xyz$.

\begin{lemma}
\label{lem:outertype}
If every graph $G_3(n)$ has geometric thickness three, then every $G_3(n)$ has an outer  drawing in which all tripletons have the same type.
\end{lemma}

\begin{proof}
By the previous lemma, we can assume that we have an outer drawing with sufficiently many vertices, to which we can apply Lemma~\ref{lem:Ram} with $c=6$ and $e=3$.
\end{proof}

We now successively analyze each type of outer drawing and show that, for sufficiently large $n$, each type must lead to a crossing.

\begin{figure}[t]
\centering
\includegraphics[width=2.25in]{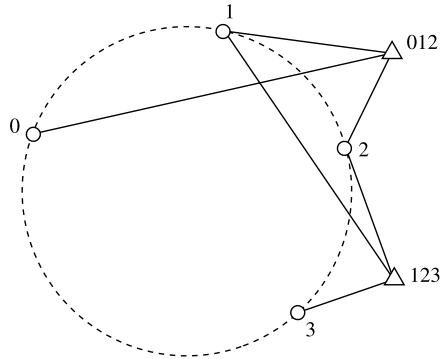}
\caption{Forced crossing in an outer drawing of type $201$, from Lemma~\ref{lem:201}.}
\label{fig:cross201}
\end{figure}

\begin{lemma}
\label{lem:201}
For $n\ge 4$, $G_3(n)$ has no outer drawing of type $201$ or $120$.
\end{lemma}

\begin{proof}
In a drawing of type $201$, each low edge of a tripleton must cross the convex curve $C$ between the middle and high sigletons adjacent to the tripleton.  Therefore, in a drawing of $G_3(4)$,
the low edge from tripleton $012$ would cross the low edge from tripleton $123$, as shown in Figure~\ref{fig:cross201}.
Type $120$ is equivalent to type $201$ under mirror reversal of the drawing and reversal of the clockwise ordering of the singletons.
\end{proof}

\begin{figure}[t]
\centering
\includegraphics[width=2in]{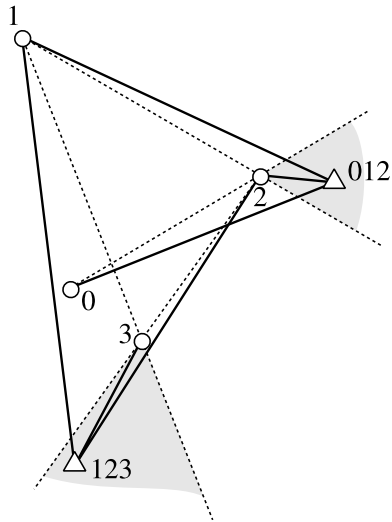}
\caption{Sharp angles in an outer drawing of type $021$, from Lemma~\ref{lem:sharp021}.}
\label{fig:sharp021}
\end{figure}

\begin{lemma}
\label{lem:sharp021}
In any outer drawing of $G_3(4)$ of type $021$, let $\theta_i$ denote the angle
at singleton $i$ on the convex hull of the singletons.  Then $\theta_1+\theta_2<180^\circ$.
\end{lemma}

\begin{proof}
In order to be of type $021$, tripleton $012$ must lie in a wedge outside $C$ bounded by lines $02$ and $12$ (the shaded area on the upper right of Figure~\ref{fig:sharp021}).
The low edge $e_0$ from $0$ to $012$ crosses line $12$ on the boundary of this wedge.
Similarly, the low edge $e_1$ from $1$ to tripleton $123$ crosses line $23$.
In order for $e_1$ to reach this line without crossing $e_0$,
tripleton $123$ must be clockwise from singleton $0$ as viewed from singleton $1$, so
angle $123$-$1$-$2$ must be greater than $\theta_1$.
But in order for the crossing of $e_1$ with line $23$ to be on the boundary of the wedge containing $123$ (shown as the lower shaded area in the figure),
angle $123$-$1$-$2$ must be less than $180^\circ - \theta_2$.
Combining these two inequalities gives the result.
\end{proof}

\begin{lemma}
\label{lem:021}
For $n\ge 6$, $G_3(n)$ has no outer drawing of type $021$ or $102$.
\end{lemma}

\begin{proof}
By applying Lemma~\ref{lem:sharp021} twice, $\theta_1+\theta_2+\theta_3+\theta_4<360^\circ$, but this is not possible in a convex polygon.
Type $102$ is equivalent to type $021$ under mirror reversal.
\end{proof}

\begin{lemma}
\label{lem:mid012}
In an outer drawing of $G_3(n)$ of type $012$, the middle edges adjacent to singletons $2$, $3$, $\ldots$, $n-3$ all pass across line segment $0$-$(n-1)$.
\end{lemma}

\begin{figure}[t]
\centering
\includegraphics[width=4.5in]{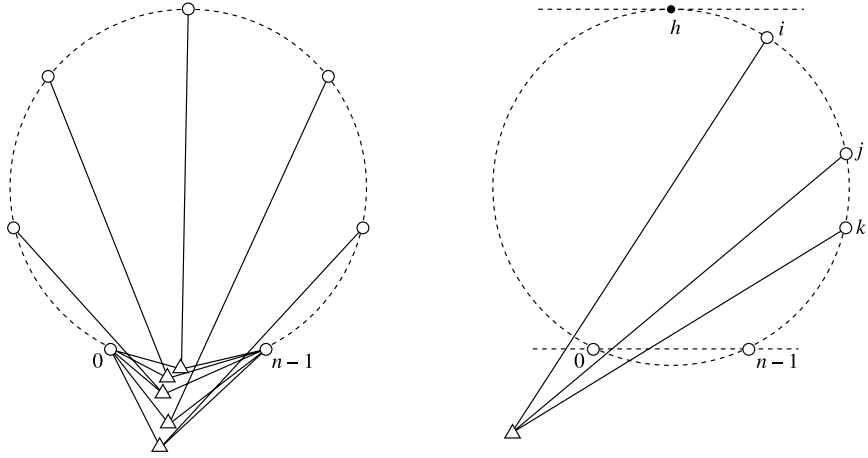}
\caption{Analysis of outer drawings of type $012$.  Left, the middle edges of the tripletons $0i(n-1)$ form tracks that force the remaining middle edges to cross segment $0$-$(n-1)$ (Lemma~\ref{lem:mid012}).  Right, low edges from clockwise of the tangent point $h$ cross $C$ counterclockwise of the tangent point (Lemma~\ref{lem:low012}).}
\label{fig:mid012}
\end{figure}

\begin{proof}
By the assumption that the drawing is of type $012$, the middle edges of the tripletons $0i(n-1)$ all pass across this segment (Figure~\ref{fig:mid012}, left).  In any other singleton $abc$, the middle edge must pass from $b$ across line segment $ac$ within the convex curve $C$.  Unless $b=1$ or $b=n-2$, this crossing with $ac$ must occur within a strip bounded by two of the middle edges of tripletons $0i(n-1)$, and the only way for the middle edge from $b$ to exit $C$ is at the end of this strip, where it meets segment $0$-$(n-1)$.
\end{proof}

\begin{lemma}
\label{lem:low012}
Suppose we have an outer drawing of $G_3(n)$ of type $012$, in which all middle edges cross segment $0$-$(n-1)$. Draw a tangent line to convex curve $C$ parallel to line $0$-$(n-1)$, on the opposite side of $C$ from that line, and let $h$ be a point where the tangent line touches $C$.  Then, for every tripleton $ijk$ where $i$ is clockwise from $h$,
the low edge from $i$ to $ijk$ crosses $C$ either at a point counterclockwise from $h$ or between singletons $0$ and $n-1$.
\end{lemma}

\begin{proof}
By the assumption of the drawing type, all low edges must cross line $0$-$n-1$ to the left of the crossing point of the middle edge, which is between points $0$ and $n-1$.  If the low edge's crossing is not also between points $0$ and $n-1$, it must be to the left of point $0$.  But, within the strip between line $0$-$(n-1)$ and the parallel tangent line, the portion of $C$ between $0$ and $h$ separates the rest of $C$ (containing $i$) from the points left of point $0$ (as shown in Figure~\ref{fig:mid012}, right), so the low edge most cross $C$ before it reaches its crossing with line $0$-$(n-1)$.
\end{proof}

Symmetrically, the high edges from counterclockwise of $h$ must cross $C$ clockwise of $h$.

\begin{lemma}
\label{lem:012}
For sufficiently large $n$, $G_3(n)$ has no outer drawing of type $012$.
\end{lemma}

\begin{proof}
Start with an outer drawing of type $012$, and remove singletons $1$ and $n-2$ if necessary so that (by Lemma~\ref{lem:mid012}) all middle edges cross segment $0$-$(n-1)$.
Form the tangent point $h$ described in Lemma~\ref{lem:low012}, and assume without loss of generality that at least half of the remaining singletons lie clockwise of $h$ on convex curve $C$. (The case that at least half the singletons lie counterclockwise of $h$ is symmetric.)
Remove from the drawing all points counterclockwise of $h$.
The result is an outer drawing of type $012$ of a graph $G_3(m)$, with at least $m\ge (n-2)/2$ singletons,
in which all low and middle edges cross segment $0$-$(m-1)$.
We can now apply the same reasoning as in Lemma~\ref{lem:noinner}:
by choosing an appropriate set of tripletons, we can find a set of low and middle edges forming a grid that contains a complete $2\times 2$ subgrid, and use the subgrid to constrain the locations of two more tripletons, forcing two high edges to cross.
\end{proof}

\begin{lemma}
\label{lem:mid210}
In an outer drawing of $G_3(n)$ of type $210$, each middle edge is completely exterior to convex curve $C$.
\end{lemma}

\begin{proof}
Suppose to the contrary that a middle edge for tripleton $ijk$ crosses $C$ at point $x$.
If $x$ is clockwise of $j$, an edge from $ijk$ to $C$ clockwise of the middle edge can only reach singletons in the range from $j$ to $x$, contradicting the requirement (from the assumption of the drawing type) that the low edge $i$-$ijk$ must be clockwise of the middle edge.  Symmetrically, if $x$ is counterclockwise of $j$, an edge counterclockwise of the middle edge can only reach singletons in the range from $x$ to $j$, contradicting the requirement that the high edge must be clockwise of the middle edge.
\end{proof}

The same reasoning as above also proves the following lemma:

\begin{lemma}
\label{lem:low210}
In an outer drawing of $G_3(n)$ of type $210$, suppose that the low edge of tripleton $ijk$ crosses the convex curve $C$.  Then the crossing point is clockwise of singleton $i$.
\end{lemma}

\begin{figure}[t]
\centering
\includegraphics[width=2in]{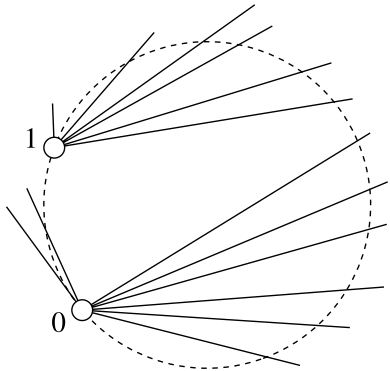}
\caption{In an outer drawing of type $210$, low edges from singletons $0$ and $1$ cross  disjoint ranges of the convex curve $C$ (Lemma~\ref{lem:nocross210}).}
\label{fig:low210}
\end{figure}

\begin{lemma}
\label{lem:nocross210}
Suppose that $G_3(n)$ has a convex drawing of type $210$, and let $n=2^k-2$  Then we can select a subset $S$ of the singletons, with $\abs{S}=k$, such that no member of $S$ is incident to a low edge that crosses the convex hull of $S$.
\end{lemma}

\begin{proof}
We use induction on $k$.  As a base case, for $k=n=2$, $G_3(2)$ has no edges.  Otherwise, consider the sets of low edges incident to singletons $0$ and $1$ that cross convex curve $C$.  Since these two sets of low edges do not cross each other (Figure~\ref{fig:low210}), we can partition the portion of $C$ clockwise of~$1$ into two parts, where the clockwise part of~$C$ is crossed only by low edges from~$0$ and the counterclockwise part of~$C$ is crossed only by low edges from~$1$.  One of these two parts contains at least $(n-2)/2$ singletons. If the clockwise part contains many singletons, we apply the induction hypothesis to these singletons, and form~$S$ by adding singleton~$1$ to the resulting set.  On the other hand, if the counterclockwise part contains many singletons, we add singleton~$0$ to the set formed by applying the induction hypothesis to the singletons in this part.

In either case, we get a set $S$ of $(k-1)+1=k$ singletons.  The low edges from $0$ or $1$ (whichever was included in the set) can not cross the hull, because of how we chose which part of $C$ to apply the induction hypothesis to.  The low edges from the remaining members of $S$ can not cross the hull between $0$ or $1$ and the rest of $S$, by Lemma~\ref{lem:low210}, and they can not cross the hull elsewhere by induction.  Therefore $S$ fulfills the conditions of the lemma.
\end{proof}

\begin{lemma}
\label{lem:210}
Not every graph $G_3(n)$ has an outer drawing of type $210$.
\end{lemma}

\begin{proof}
By Lemma~\ref{lem:nocross210}, if $G_3(n)$ has such a drawing, then we can find
a set of $k$ singletons, and a convex curve $C$ (the hull of this set),
forming an outer drawing of $G_3(k)$ in which
all low and middle edges avoid the interior of convex curve $C$;
here $k=\lfloor \log_2(n+2)\rfloor$.
We can now apply the same reasoning as in Lemmas \ref{lem:noinner} and~\ref{lem:012}:
by choosing an appropriate set of tripletons, we can find a set of low and middle edges forming a grid that contains a complete $2\times 2$ subgrid, and use the subgrid to constrain the locations of two more tripletons, forcing two high edges to cross.
\end{proof}

\section{Completion of the Proof}

\begin{lemma}
\label{lem:notgt3}
Not every graph $G_3(n)$ has geometric thickness at most three.
\end{lemma}

\begin{proof}
This follows from Lemma~\ref{lem:outertype} (showing that we can restrict our attention to outer drawings of fixed type) and Lemmas \ref{lem:201}, \ref{lem:021}, \ref{lem:012}, and \ref{lem:210}
(showing that each fixed type leads to a forced crossing).
\end{proof}

\begin{figure}[t]
\centering
\includegraphics[width=2.5in]{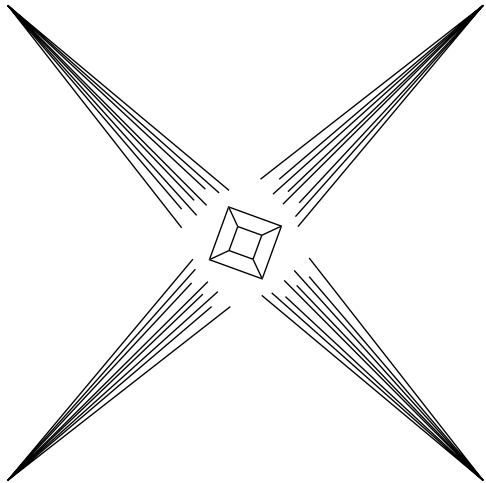}
\caption{Schematic view of drawing of  $G_3(8)$ with geometric thickness three.}
\label{fig:g38}
\end{figure}

The largest $n$ for which $G_3(n)$ has no thickness-three drawing must be at least nine, as $G_3(8)$ has geometric thickness at most three: draw four singleton subsets in a large square, with the remaining four singletons and all tripletons clustered near the center of the square.  The edges incident to the inner four singletons form a planar graph, consisting of a cube $G_3(4)$, 24 two-edge paths across the faces of the cube, and 24 tripletons with only a single edge to one of the inner singletons.  By choosing an appropriate F\'ary embedding of this planar graph, we can draw it in a single layer, and use an additional layer for the edges incident to each pair of opposite outer singletons.  Figure~\ref{fig:g38} shows this drawing in a schematic view, with the outer four vertices and inner cube visible but omitting the remaining tripletons.

\section{Conclusions}

We have shown that thickness and geometric thickness are not asymptotically equivalent.
Some interesting open questions remain:

\begin{enumerate}
\item
How does $\bar\theta(G_3(n))$ grow as a function of $n$?
The best upper bound we have been able to obtain is that
$\bar\theta(G_3(n))\le \lceil(n-2)/2\rceil$, by generalizing the three-layer drawing of $G_3(8)$ outlined
in the previous section, but this is far from the lower bound given by
Theorem~\ref{thm:sepgeom}.

\item
What can be said about two-layer drawings?
It is known that there exist graphs, for instance $K_{6,8}$, for which $\theta(G)=2$ and $\bar\theta(G)>2$~\cite{Man-MPCPS-83,DilEppHir-JGAA-00}.  Can $\bar\theta(G)$ be unbounded when $\theta(G)=2$?

\item
Four-regular graphs have thickness at most two.  Do random four-regular graphs have bounded expected geometric thickness? More generally, do bounded degree graphs have bounded geometric thickness?

\item
The two-layer drawing technique for $G_2(n)$
generalizes to many other {\em $2$-inductive} graphs: that is, graphs that can be reduced to the empty graph by repeatedly removing vertices with degree at most two.
Is there a $2$-inductive graph with geometric thickness greater than two?
Can the geometric thickness of such graphs be unbounded?

\item
Due to our use of Ramsey theory, our lower bounds on geometric thickness grow only very slowly as a function of the number of vertices.  Can we find good upper bounds for geometric thickness as a combination of the graph-theoretic thickness and a slowly growing function of the graph size?

\item What is the complexity of computing $\bar\theta(G)$ for a given graph $G$~\cite{DilEppHir-JGAA-00}?

\item Thickness is NP-hard~\cite{Man-MPCPS-83} but not difficult to approximate to within a constant factor; e.g., it is within a factor of three of the graph's arboricity.  Our new results imply that these approximations do not directly extend to geometric thickness.  Is there an efficient algorithm for layered drawing of graphs without bends, using a number of layers within a constant factor of optimal?
\end{enumerate}

\raggedright
\bibliographystyle{amsplain}
\bibliography{thickness}
 \end{document}